\documentclass{article}
\usepackage{amsmath}
\usepackage{amsfonts}
\usepackage{amsthm}

\newfont{\rus}{wncyr10 scaled 1200}
\newfont{\rusb}{wncyb10 scaled 1200}

\newtheorem{theorem}{Theorem}

\newtheorem{lemma}[theorem]{Lemma}

\newcommand{\lin}{{\rm lin}\hskip0.02cm}

\newcommand{\ep}{\varepsilon}

\begin{document}
\title{\LARGE{\bf Weak$^*$ closures and derived sets in dual Banach spaces}}

\author{M.\,I.~Ostrovskii\\
Department of Mathematics and Computer Science\\
St. John's University\\
8000 Utopia Parkway\\
Queens, NY 11439\\
USA\\
e-mail: {\tt ostrovsm@stjohns.edu}}

\date{\today}
\maketitle

\begin{large}

\hfill{\it Dedicated to the memory of Vincenzo Bruno Moscatelli}
\bigskip

\noindent{\bf Abstract:} The main results of the paper: {\bf (1)}
The dual Banach space $X^*$ contains a linear subspace $A\subset
X^*$ such that the set $A^{(1)}$ of all limits of weak$^*$
convergent bounded nets in $A$ is a proper norm-dense subset of
$X^*$ if and only if $X$ is a non-quasi-reflexive Banach space
containing an infinite-dimensional subspace with separable dual.
{\bf (2)} Let $X$ be a non-reflexive Banach space. Then there
exists a convex subset $A\subset X^*$ such that $A^{(1)}\neq
{\overline{A}\,}^*$ (the latter denotes the weak$^*$ closure of
$A$). {\bf (3)} Let $X$ be a quasi-reflexive Banach space and
$A\subset X^*$ be an absolutely convex subset. Then
$A^{(1)}={\overline{A}\,}^*$.
\medskip

\medskip

\medskip

Let $A$ be a subset of a dual Banach space $X^*$, we denote the
weak$^*$ closure of $A$ by ${\overline{A}\,}^*$. The {\it weak$^*$
derived set} of $A$ is defined as
\[A^{(1)} = \bigcup_{n=1}^\infty {\overline{A\cap
nB_{X^*}}\,}^{*},\] where $B_{X^*}$ is the unit ball of $X^*$,
that is, $A^{(1)}$ is the set of all limits of weak$^*$ convergent
bounded nets in $A$. If $X$ is separable, $A^{(1)}$ coincides with
the set of all limits of weak$^*$ convergent sequences from $A$,
called the {\it weak$^*$ sequential closure}. The strong closure
of a set $A$ in a Banach space is denoted $\overline{A}$. A subset
$A\subset X^*$ is called {\it total} if for every $0\ne x\in X$
there exists $f\in A$ such that $f(x)\ne 0$. A subset $A\subset
X^*$ is called {\it norming} if there is $c>0$ such that for every
$0\ne x\in X$ there exists $f\in A$ satisfying $||f||=1$ and
$f(x)\ge c||x||$.
\medskip

The study of weak$^*$ derived sets was initiated by Banach and
continued by many authors, see
\cite{Maz30,Ban32,McG68,Sar68,God78,Mos87,Ost87}, and references
therein. Weak$^*$ derived sets and their relations with weak$^*$
closures found applications in many areas: the structure theory of
Fr\'echet spaces (see
\cite{Alb05,BDH86,DM87,MM89,MM92,Mos90,Ost98}), Borel and Baire
classification of linear operators, including the theory of
ill-posed problems (\cite{PP80,Pli97,Raj04,Sai76}), Harmonic
Analysis (\cite{KL87,Lyo88,McG68,Pia54}), theory of biorthogonal
systems (\cite{HMVZ07,Pli86}; I have to mention that the
historical information on weak$^*$ sequential closures in
\cite{HMVZ07} is inaccurate). The survey \cite{Ost01} contains a
historical account and an up-to-date-in-2000 information on
weak$^*$ sequential closures.
\medskip

Recently  derived sets were used in the study of extension
problems for holomorphic functions on dual Banach spaces
\cite{GKM10}.
\medskip

The main purpose of the present paper is to answer the following
two questions asked in \cite{GKM10}:\medskip

\noindent{\bf 1.} \cite[Question 6.3(a)]{GKM10}. Let $X$ be a
quasi-reflexive Banach space.  Is ${\overline{A}\,}^{*} = A^{(1)}$
for each (absolutely) convex set $A\subset X^*$?\medskip

\noindent{\bf 2.} \cite[Question 6.5]{GKM10}. For which Banach
spaces $X$ there is a linear subspace $A\subset X^*$ such that
$A^{(1)}$ is a proper norm-dense subset of $X^*$? Is it true
whenever $X$ is not quasi-reflexive?
\medskip

The main results of the paper:

\begin{theorem}\label{T:NonQuasi} The dual Banach space $X^*$ contains a linear subspace $A\subset X^*$ such that
$A^{(1)}$ is a proper norm-dense subset of $X^*$ if and only if
$X$ is a non-quasi-reflexive Banach space containing an
infinite-dimensional subspace with separable dual.
\end{theorem}

\begin{theorem}\label{T:NonRefl} Let $X$ be a non-reflexive Banach space. Then there
exists a convex subset $A\subset X^*$ such that $A^{(1)}\neq
{\overline{A}\,}^*$.
\end{theorem}

\begin{theorem}\label{T:QuasiRefl} Let $X$ be a quasi-reflexive Banach space and $A\subset X^*$
be an absolutely convex subset. Then $A^{(1)}={\overline{A}\,}^*$.
\end{theorem}

Some parts of Theorems \ref{T:NonQuasi} and \ref{T:NonRefl} are
proved for separable spaces with basic sequences of special kinds
first, and then are extended to the general case.\medskip

To describe the way in which results are extended from subspaces
we need some more notation. Let $Z$ be a subspace in a Banach
space $X$ and $E:Z\to X$ be the natural isometric embedding. Then
$E^*:X^*\to Z^*$ is a quotient mapping which maps each functional
in $X^*$ onto its restriction to $Z$. Let $A$ be a subset of
$Z^*$. It is clear that $D=(E^*)^{-1}(A)$ is the set of all
extensions of all functionals in $A$ to the space $X$.

\begin{lemma}\label{L:DerExt}
\begin{equation}\label{E:DerExt}
D^{(1)}=(E^*)^{-1}(A^{(1)}),
\end{equation}
where the derived set $D^{(1)}$ is taken in $X^*$ and the derived
set $A^{(1)}$ - in $Z^*$.
\end{lemma}

\begin{proof} The inclusion $D^{(1)}\subset(E^*)^{-1}(A^{(1)})$ follows from the inclusion
$E^*(D^{(1)})\subset A^{(1)}$, which in turn follows from the
weak$^*$ continuity of $E^*$.
\medskip

To  prove  the  inverse inclusion it suffices to show that for
every bounded net $\{f_\nu\}\subset Z^{*}$ with $w^*-\lim_\nu
f_\nu=f$ and every $g\in (E^{*})^{-1}(\{f\})$ there exist
$g_{\nu}\in (E^{*})^{-1}(\{f_\nu\})$ such that some subnet of
$\{g_\nu\}$ is bounded and weak$^*$ convergent to $g$. Let
$h_{\nu}$ be such that $h_{\nu}\in (E^{*})^{-1}(\{f_{\nu}\})$ and
$||h_{\nu}||= ||f_{\nu}||$ (Hahn-Banach extensions). Then
$\{h_{\nu}\}_\nu$ is a bounded net in $X^{*}$. Hence it has a
weak$^*$ convergent subnet, let $h$ be its limit. Then
$g-h\in(E^{*})^{-1}(\{0\})$, therefore $g_\nu=h_\nu+g-h$ is a
desired net.
\end{proof}

\begin{proof}[Proof of Theorem \ref{T:NonQuasi}] First we suppose
that $X$ is such that $X^*$ contains a subspace $A$ for which
$A^{(1)}$ is a proper norm-dense subset in $X^*$. \medskip

The space $X$ cannot be quasi-reflexive because the norm-density
of $A^{(1)}$ in $X^*$ implies that $A$ is total, and the condition
$A^{(1)}\ne X^*$ implies that $A$ is not norming \cite{Dix48}, and
total non-norming subspaces do not exist in duals of
quasi-reflexive spaces (\cite{Pet64}, \cite{Sin63}).
\medskip

To show that $X$ contains an infinite-dimensional subspace with
separable dual, assume the contrary, that is, all
infinite-dimensional subspaces of $X$ have non-separable duals.
\medskip

Define the Banach space $X_A$ as the completion of $X$ with
respect to the norm $||x||_A=\sup\{|f(x)|:~f\in A,~||f||=1\}$.
Since the subspace $A$ is non-norming, the natural mapping $N:X\to
X_A$ is not an isomorphism. Since the subspace $A$ is total, the
mapping $N$ is injective.\medskip

Using the standard argument \cite[Proposition 2.c.4]{LT77} we find
a separable infinite-dimensional subspace $Z\subset X$ such that
the restriction $N|_Z$ is a compact operator. By duality
\cite[VI.5.2]{DS58}, this implies that $R=(N|_Z)^*(B_{X_A^*})$ is
a norm-compact subset of $Z^*$. Observe that $A \cap B_{X^*}$ is
embedded in a natural way into $B_{X_A^*}$. Therefore $E^*(A \cap
B_{X^*})\subset R$. Therefore $E^*$ maps each weak$^*$ convergent
net in $A\cap B_{X^*}$ onto a strongly convergent net in $Z^*$,
therefore $E^*(A^{(1)})$ is contained in the linear span of $R$,
which is a separable subspace of $Z^*$. Since by our assumption
$Z^*$ is non-separable, the subspace $E^*(A^{(1)})$ is not dense
in $Z^*$. Hence $A^{(1)}$ is not dense in $X^*$, this
contradiction completes the first part of the proof.
\medskip

Now we prove the converse. Assume that $X$ is a
non-quasi-reflexive Banach space containing an
infinite-dimensional subspace with separable dual. We use
terminology of \cite{LT77}. The following result is proved using
the techniques of \cite{DJ73}. We use the notation
$n_k=\frac{k(k+1)}2$ for $k=0,1,2,\dots$.

\begin{lemma}\label{L:GenDJ} Let $X$ be a non-quasi-reflexive Banach spaces
containing an infinite-dimensio\-nal subspace with separable dual.
Then there exists a minimal system
\[\{u_i\}_{i=0}^\infty\cup \{x_i\}_{i=0}^\infty\]
in $X$ satisfying the conditions:

\begin{itemize}

\item[{\bf (1)}] The system $\{u_i\}_{i=0}^\infty\cup
\{x_i\}_{i=0}^\infty$ and its biorthogonal functionals
$\{u^*_i\}_{i=0}^\infty\cup \{x^*_i\}_{i=0}^\infty$ are uniformly
bounded.

\item[{\bf (2)}] The sequence $\{u_i\}_{i=0}^\infty$ spans a
subspace $U$ with separable dual $U^*$, and the restrictions of
the biorthogonal sequence $\{u_i^*\}_{i=0}^\infty$ to $U^*$ span
$U^*$.

\item[{\bf (3)}] The set $\left\{\sum_{p=j}^kx_{n_p+j}:~ 0\le j\le
k\le\infty\right\}$ is bounded.
\end{itemize}
\end{lemma}

\begin{proof} This proof is a modification of the proof of
Proposition 1 in \cite[pp.~360--362]{DJ73}. For this reason we
mostly follow the terminology and notation of \cite{DJ73} and the
reader is expected to consult \cite{DJ73} if more details are
needed (making this proof readable independently from \cite{DJ73}
would lead to too much copying from \cite{DJ73}).\medskip

For the same reason as in \cite[Proposition 1]{DJ73} we may assume
that $X$ is a separable non-quasi-reflexive Banach space
containing an infinite-dimensional subspace $Y$ with separable
dual. Let sequence $\{s_i\}_{i=0}^\infty\subset X^*$ be such that
its restrictions to $Y^*$ are dense in $Y^*$. By \cite[Theorem
1]{DJ73}, there is a weak$^*$ null sequence $\{y_n\}\subset X^*$,
a bounded sequence $\{f_n\}_{n=0}^\infty\subset X^{**}$, and a
partition $I_n$ of the integers into pairwise disjoint infinite
subsets such that
\[f_k(y_n)=\begin{cases} 1 &\hbox{ if }n\in I_k\\
0 &\hbox{ if }n\notin I_k.\end{cases}
\]
Let $\lambda>||f_n||$ for all $n$ and choose $0<\ep_i<1$ for all
$i$ with $\prod_i(1+\ep_i)<\infty$.\medskip

We are going to use induction to show that for each $p$ we can
find $\{x_k|~0\le k\le n_p + p\}\subset X$, $\{f_i'|~0\le i\le
p\}\subset X^{**}$, $\{u_i|~0\le i\le p\}\subset Y$, and
finite-dimensional subspaces $G_0\subset G_1\subset\dots\subset
G_p$ of $X^*$ such that the following conditions are satisfied:

\begin{enumerate}

\item $G_i$ is $(1 +\ep_i)^2$ norming over the linear span of
$\{x_k|~ 0\le k\le n_i+ i\}\cup\{u_j:0\le j\le i\}$ for $0\le i\le
p$.

\item $s_i\in G_i$ for each $i=0,\dots,p$.

\item $u_i\in Y\cap (G_{i-1})_\bot$ for $1\le i\le p$,
$||u_i||=1$.

\item $\{x_{n_{i}+j}|~0\le j\le i\}\subset(G_{i-1})_\bot$ for
$1\le i\le p$.

\item $\displaystyle{\left(\left.\left\|\sum_{i=j}^k
x_{n_i+j}\right\|\ \right|~j\le k\le p\right)}$ is bounded by
$6\lambda$, for $0\le j\le p$.

\item $g\left(\sum^p_{i=j} x_{n_i+j}\right)=f'_j(g)$ for $g\in
G_p$, $0\le j\le p$.

\item There is a constant $C$ depending only on $\sup_n||y_n||$
such that for each $j$, $0\le j\le p$ there are functionals
$\{\varphi_{n_j+i}|~0\le i\le j\}\subset X^*$ of norm $\le C$ such
that the system $\{x_{n_j+i},\varphi_{n_j+i}|~0\le i\le j\}$ is
biorthogonal.

\item For each $j$, $0\le j\le p$ there is a functional $v_j^*\in
X^*$ such that $||v_j^*||=1$, $v_j^*(u_j)=1$, and
$v_j^*(x_{n_j+i})=0$ for $0\le i\le j$.

\item $||f_i'||\le 3\lambda$ for $0\le  i\le p$.

\item There exist infinite sets $I_k'\subset I_k$, $k = 0,
1,\dots$, so that for each $i, 0\le i\le p$, $f_i'$ agrees with
$f_i$ on $[y_n| n\in\cup I_k']$. The sets $I_k'$ depend on $i$,
although we do not reflect this dependence in our notation.
\end{enumerate}

For the first step, let $U_0$ be a $2$-dimensional subspace of $Y$
and $G_0$ be a finite-dimensional subspace of $X^*$ which $(1
+\ep_0)$-norms $[\{f_0\}\cup U_0]$. By local reflexivity
\cite{JRZ71,LR69}, we pick $x_0$ in $X$ with
$||x_0||\le\min\{\lambda, (1+\ep_0)||f_0||\}$ such that $g(x_0)
=f_0(g)$ for $g$ in $G_0$. For convenience of notation later, we
rename $f_0$ by $f_0'$. By the well-known result of \cite{KKM48}
(see \cite[Lemma 2.c.8]{LT77}) there is  $u_0\in U_0$, $||u_0||=1$
such that for some $v^*_0\in X^*$ we have $||v^*_0||=1$,
$v_0(u_0)=1$ and $v_0(x_0)=0$.
\medskip

Let $(1'),\dots,(10')$ be the statements above for $p+1$. By
\cite[Lemma 1]{DJ73}, pick infinite sets $I_k''\subset I_k'$ for
all $k$ so that the natural projection onto $G_p$ from $G_p\oplus
[y_n|~ n\in\cup I_k'']$ has norm $\le 2$. Hence, there exists
$f'_{p+1}$ in $X^{**}$ with $||f'_{p+1}||<3\lambda$ so that
$f'_{p+1}(g) = 0$ for $g\in G_p$, and such that $f'_{p+1}$ agrees
with $f_{p+1}$ on $[y_n|~n\in\cup I_k'']$. This satisfies $(9')$
and $(10')$.\medskip

Since $y_n\overset{w^*}\to0$, and each $I_k''$ is infinite, there
exist, for $0\le i\le p + 1$, $q_i\in I_i''$ so that
$\sum_{k=0}^{n_p+p}|y_{q_i}(x_k)|< 1/4p$. Now we select a
$(p+2)$-dimensional subspace $U_{p+1}\subset Y\cap (G_p)_\bot$ and
a finite-dimensional subspace $G_{p+ 1}\subset X^*$ containing
$G_p\cup\{s_{p+1}\}\cup \{y_{q_i}|~ 0 \le i \le p + 1\}$ and such
that $G_{p+1}$ is $(1 + \ep_{p+1})$-norming over the linear span
$H$ of $\{x_k\}_{k=0}^{n_p+p}\cup \{f_i'\}_{i=0}^{p+
1}\cup\{u_i\}_{i=0}^p\cup U_{p+1}$ in $X^{**}$. This definition of
$G_{p+1}$ implies that $(2')$ is satisfied. By the principle of
local reflexivity \cite{JRZ71,LR69}, there is an operator $T:H\to
X$ such that $T$ is the identity on
$\{x_k\}_{k=0}^{n_p+p}\cup\{u_i\}_{i=0}^p\cup U_{p+1}$, $T$ is an
$(1 + \ep_{p+1})$-isometry and $g(Tf)=f(g)$ for $f\in H$, $g\in
G_{p+1}$. Define $x_{n_{p+1}},\dots,x_{n_{p+1}+p+1}$ by
$x_{n_{p+1}+j}= Tf_j'-\sum_{i=j}^px_{n_i+j}$ for $0\le j\le p$ and
$x_{n_{p+1}+p+1}=Tf'_{p+1}$.\medskip

Thus $Tf_j'=\sum_{i=j}^{p+1}x_{n_i+j}$ for $0\le j\le p + 1$, so
that $(5')$ and $(6')$ hold. Since $G_p\subset G_{p+ 1}$ and
$f'_{p+1}\in G_p^\perp$, using (6) we get $(4')$. Now, for $0\le
i\le p + 1$, $0\le j\le p + 1$, one again has from local
reflexivity that
$y_{q_i}(x_{n_{p+1}+j})=f_j'(y_{q_i})-\sum_{k=j}^p
y_{q_i}(x_{n_k+j})$, so that $y_{q_i}(x_{n_{p+1}+i})\ge 1 -
1/4p\ge 3/4$ and $|y_{q_i}(x_{n_{p+1}+j})| < 1/4p$ when $i\ne j$.
It is easy to derive from these inequalities that the Hahn-Banach
extensions of the functionals defined by
$\varphi_{n_p+i}(x_{n_p+j})=\delta_{ij}$, $i,j=0,1, ...,p + 1$,
satisfy $||\varphi_{n_p+i}||\le C$ where $C$ depends on
$\sup_n||y_n||$ only, so $(7\,')$ is satisfied.
\medskip

Now we use \cite[Lemma 2.c.8]{LT77} and pick $u_{p+1}\in U_{p+1}$
such that $||u_{p+1}||=1$ such that for some $v^*_{p+1}\in X^*$ we
have $||v^*_{p+1}||=1$, $v_{p+1}(u_{p+1})=1$ and
$v_{p+1}(x_{n_{p+1}+i})=0$ for $i=1,\dots,p+1$. It is clear that
$(3')$ and $(8')$ are satisfied.\medskip

Since $G_{p+1}$ is $(1 + \ep_{p+1})$-norming over $H$, local
reflexivity guarantees that $G_{p+1}$ is $(1
+\ep_{p+1})^2$-norming over the linear span of $\{x_k|~ 0\le k\le
n_{p+1}+ p+1\}\cup\{u_i:0\le i\le p+1\}$ so that $(1')$ holds.
This completes the construction of $\{x_n\}$ and
$\{u_n\}$.\medskip

The conditions (1), (3), and (4) and the choice of $\{\ep_i\}$
imply that the sequence
\[[x_0,u_0], [x_1,x_2,u_1], \dots, [x_{n_p},\dots,x_{n_p+p},u_p], \dots
\]
of subspaces forms a finite-dimensional decomposition of the
closed linear span of
\[ \{u_i\}_{i=0}^\infty\cup \{x_i\}_{i=0}^\infty
\]

Now we check that conditions {\bf (1--3)} of Lemma \ref{L:GenDJ}
are satisfied.\medskip

The sequences $\{x_i\}$ and $\{u_i\}$ are bounded by construction,
for $\{x_i\}$ we use also (5). The biorthogonal functionals of the
system $\{u_i\}_{i=0}^\infty\cup \{x_i\}_{i=0}^\infty$ are bounded
because of the finite decomposition property and the conditions
(7) and (8). It remains to check the condition {\bf (2)}.\medskip

Since $\{u_i\}_{i=0}^\infty$ is a basis in its closed linear span,
it suffices to show that this basis is shrinking, that is, that
for each $u^*\in U^*$ and $\delta>0$ there is $n\in\mathbb{N}$
such that $||u^*|_{[u_i]_{i=n}^\infty}||<\delta$ (see
\cite[Proposition 1.b.1]{LT77}). Let $\tilde u$ be a
norm-preserving extension of $u^*$ to $Y$. By density there exists
$n\in\mathbb{N}$ such that $||s_{n-1}|_Y-\tilde u||<\delta$. The
conditions (2) and (3) above imply that
$s_{n-1}|_{[u_i]_{i=n}^\infty}=0$. Hence
$||u^*|_{[u_i]_{i=n}^\infty}||<\delta$.
\end{proof}

We consider the subspace $W$ spanned by the system
$\{u_i\}_{i=0}^\infty\cup \{x_i\}_{i=0}^\infty$ constructed in
Lemma \ref{L:GenDJ}. Denote by $h_j$, $j\in\mathbb{N}\cup\{0\}$, a
weak$^*$-cluster point of the sequence \[{\left\{\sum_{i=j}^k
x_{n_i+j}\right\}_{k=j}^\infty}\] in $W^{**}$, and by
$\{u^*_i\}_{i=0}^\infty\cup \{x^*_i\}_{i=0}^\infty \subset W^*$
the biorthogonal functionals of $\{u_i\}_{i=0}^\infty\cup
\{x_i\}_{i=0}^\infty$.
\medskip

It is easy to see that Lemma \ref{L:DerExt} and the Hahn-Banach
theorem imply that it suffices to find a subspace $A\subset W^*$
such that $A^{(1)}\ne W^*$ and $A^{(1)}$ is norm-dense in $W^*$.
\medskip

To construct such $A$ we pick a sequence $\{a_i\}_{i=1}^\infty$ of
real numbers satisfying $a_i>0$ and $\sum_{i=1}^\infty
a_i<\infty$, and let
$K=\{u_i+a_ih_i:~i\in\mathbb{N}\cup\{0\}\}\subset W^{**}$ and
$A=K_{\bot}\subset W^*$.\medskip

We claim that
\medskip

(A) $u_i^*\in A^{(1)}$ for all $i$.
\medskip

In fact, $u_i^*$ is a weak$^*$ limit of
$u^*_i-\frac{1}{a_i}x^*_{n_j+i}$ as $j\to\infty$ and
$u^*_i-\frac{1}{a_i}x^*_{n_j+i}\in K_{\bot}$ for $j\ge i$.
\medskip

(B) If $y^*\in U^{\perp}\subset W^*$, then
\[y^*-\sum_ia_ih_i(y^*)u_i^*\in A.\]
This immediately follows from the condition $h_i(u_j^*)=0$. The
series is norm-convergent because $\{h_j\}$ and $\{u_i^*\}$ are
bounded sequences and $\sum_{i=1}^\infty a_i<\infty$. Therefore
$y^*\in\overline{A^{(1)}}$.\medskip

We show that conditions (A) and (B) imply that
$\overline{A^{(1)}}=W^*$. In fact, let $z^*\in W^*$, $\ep>0$. By
Lemma \ref{L:GenDJ}{\bf (2)}, the restriction of $z^*$ to $U$ can
be $\ep$-approximated by a finite linear combination of
restrictions of $u_i^*$ to $U$. Therefore there exists a vector
$s$ in $U^*$ such that $||s||<\ep$ and $z^*|_U-s$ is a finite
linear combination of $\{u_i^*|_U\}$. Let $s^*$ be a Hahn-Banach
extension of $s$ to $W$, so $||s^*||<\ep$ and only finitely many
of the numbers $\{(z^*-s^*)(u_i)\}_{i=1}^\infty$ are non-zero.
Subtracting from $z^*-s^*$ the corresponding finite linear
combination of $\{u^*_i\}$ we get a vector from $U^\perp$. Thus
every vector $z^*\in W^*$ can be arbitrarily well approximated by
vectors of $\lin(\{u_i^*\}\cup U^\perp)$, hence
$W^*=\overline{A^{(1)}}$.
\medskip

It remains to prove that $A^{(1)}\ne W^*$. Let $z^*\in A\cap
B_{W^*}$. Then $|z^*(u_i)|=|-a_ih_i(z^*)|\le a_iC$, where
$C=\sup_i||h_i||$. It is easy to see that the conditions on
$\{u_i\}$ and $\{a_i\}$ imply that the set $T=\{u^*\in U^*:
|u^*(u_i)|\le a_iC~\forall i\in\mathbb{N}\cup\{0\}\}$ is
norm-compact. The inequality $|z^*(u_i)|\le a_iC$ implies that
$E^*(A\cap B_{W^*})\subset T$, where $E$ is the natural isometric
embedding of $U$ into $W$. Since $T$ is norm-compact, the set
$E^*\left({\overline{A\cap B_{W^*}}\,}^*\right)$ is also contained
$T$. Therefore $E^*(A^{(1)})\subset\lin(T)\ne U^*$ and $A^{(1)}\ne
W^*$.
\end{proof}

\begin{proof}[Proof of Theorem \ref{T:NonRefl}] We are going to
use the following result proved in \cite{Pel62,Sin62}: If a Banach
space $X$ is non-reflexive, then it contains a normalized basic
sequence $\{z_i\}_{i=0}^\infty$ such that the sequence
$\left\{\sum_{i=1}^kz_i\right\}_{k=1}^\infty$ is bounded. Let $Z$
be the closed linear span of the sequence $\{z_i\}_{i=0}^\infty$
and $z^{**}$ be a weak$^*$-cluster point of the sequence
$\left\{\sum_{i=1}^kz_i\right\}_{k=1}^\infty$ in $Z^{**}$. (We
added an extra element $z_0$ to the sequence, because it is needed
for our construction, of course, it does not affect the validity
of the result of \cite{Pel62,Sin62}.) By Lemma \ref{L:DerExt}, it
suffices to find a convex subset $A\subset Z^*$ such that
$A^{(1)}\ne{\overline{A}\,}^*$. In fact, if we have such $A$, we
let $D=(E^*)^{-1}(A)$. We have, by Lemma \ref{L:DerExt},
$D^{(1)}=(E^*)^{-1}(A^{(1)})$. Also, by the bipolar theorem
${\overline{A}\,}^*=A^{\circ\circ}$, where the first polar is in
$Z$ and the second in $Z^*$. It is easy to see that the polar
$D^{\circ}$ of $D$ in $X$ coincides with $A^\circ$. Therefore
${\overline{D}\,}^*=A^{\circ\circ}$, where the first polar is in
$Z$, and the second in $X^*$. Hence
${\overline{D}\,}^*\supset(E^*)^{-1}({\overline{A}\,}^*)$ and
$D^{(1)}\ne{\overline{D}\,}^*$.
\medskip

Let $\{z_i^*\}$ be the biorthogonal functionals of $\{z_i\}$,
$\{\alpha_i\}$ and $\{\beta_i\}$ be strictly increasing sequences
of positive real numbers satisfying $\lim_{i\to\infty}\alpha_i=1$
and $\lim_{i\to\infty}\beta_i=\infty$. We split $\mathbb{N}$ into
infinitely many infinite subsequences $\mathbb{N}_j$. Let
$A\subset Z^*$ be the convex hull of all vectors of the form
$\alpha_jz_0^*+\beta_jz_k^*$, where $k\in\mathbb{N}_j$.
\medskip

It is enough to show that the set $A^{(1)}$ is not strongly
closed. First we observe that $\alpha_jz_0^*\in A^{(1)}$. In fact,
$\alpha_jz_0^*$ is the weak$^*$ limit of the sequence
$\{\alpha_jz_0^*+\beta_jz_k^*\}_{k\in \mathbb{N}_j}$.
\medskip

It remains to show that $z_0^*\notin A^{(1)}$. Assume the
contrary. Since $Z$ is separable, there is a bounded sequence
$\{y^*_r\}_{r=1}^\infty$ of vectors in $A$ such that $z_0^*$ is a
weak$^*$ limit of $\{y^*_r\}_{r=1}^\infty$. By the definition of
$A$, vectors $y^*_r$ are finite convex combinations of the form
$y^*_r=\sum_{j,k} a_{j,k}(r)(\alpha_jz_0^*+\beta_jz_k^*)$. Since
$z_0$ is a weak$^*$ continuous functional on $Z^*$, we get
\[\lim_{r\to\infty}\sum_{j=1}^\infty a_{j,k}(r)\alpha_j=1.\]

It is clear that this implies

 \[\lim_{r\to\infty}\sum_{\substack{j\\~\alpha_j<1-\ep}}
a_{j,k}(r)=0.\]

Since $\lim_{j\to\infty}\beta_j=\infty$, this implies that for
each $M<\infty$
 \[\lim_{r\to\infty}\sum_{\substack{j\\~\beta_j>M}}
a_{j,k}(r)=1.\] Therefore
\[\begin{split}\limsup_{r\to\infty}z^{**}(y_r)&=
\limsup_{r\to\infty}\sum_{j=1}^\infty a_{j,k}(r)\beta_j\\&\ge
M\limsup_{r\to\infty}\sum_{\substack{j\\~\beta_j>M}}a_{j,k}(r)=M.\end{split}\]
Since $M$ is arbitrary, this implies that the sequence $y^*_r$ is
unbounded, and we get a contradiction.
\end{proof}

\begin{proof}[Proof of Theorem \ref{T:QuasiRefl}] Assume the contrary, let $A$ be an
absolutely convex subset of the dual of a quasi-reflexive Banach
space $X$ such that $A^{(1)}\ne {\overline{A}\,}^*$.\medskip

By \cite[Theorem V.5.7]{DS58} this implies that
$\left(A^{(1)}\right)^{(1)}\ne A^{(1)}$, so there exists a bounded
weak$^*$ convergent net $\{x_\alpha\}$ in $A^{(1)}$ such that any
nets $\{x_{\alpha,\beta}\}\subset A$ satisfying
\begin{equation}\label{E:AlphaBeta}\sup_\beta||x_{\alpha,\beta}||<\infty\quad\hbox{and}\quad w^*-\lim_\beta
x_{\alpha,\beta}=x_\alpha\end{equation} are not uniformly bounded,
that is,
\[\sup_\alpha\sup_\beta||x_{\alpha,\beta}||=\infty.\]

Since $X$ is quasi-reflexive, we have $X^{**}=X\oplus F$ where $F$
is a finite-dimensional subspace. We pick nets
$\{x_{\alpha,\beta}\}_\beta\subset X^*$ satisfying the condition
\eqref{E:AlphaBeta}. We may assume that $\beta$ in all of them
runs through the same ordered set (it can be chosen to be a subnet
of the naturally ordered set of weak$^*$ neighborhoods of $0$ in
$X^*$) and that the natural images of these nets in $F^*$ converge
strongly. Denote the corresponding limits in $F^*$ by $v_\alpha$.
First we show that $\limsup_\alpha||v_\alpha||=\infty$.\medskip

Assume the contrary, that is, $\limsup_\alpha||v_\alpha||<\infty$.
Using local reflexivity \cite{JRZ71,LR69} we find, for
sufficiently large $\alpha$, uniformly bounded nets
$\{\ell_{\alpha,\delta}\}_\delta\in X^{*}$ such that
$\ell_{\alpha,\delta}|_F=v_\alpha-(x_\alpha|_F)$ and
$\lim_\delta\ell_{\alpha,\delta}(x)=0$ for all $x\in X$.

Then the combined nets
$\{x_{\alpha,\beta}-x_\alpha-\ell_{\alpha,\delta}\}_{(\beta,\delta)}$
(where the order is defined by: $(\beta_1,\delta_1)\succ
(\beta_2,\delta_2)$ if and only if both $\beta_1\succ\beta_2$ and
$\delta_1\succ\delta_2$) are weakly null. In fact, if $x\in X$
then $\lim_\beta x_{\alpha,\beta}(x)=x_\alpha(x)$ and
$\lim_\delta\ell_{\alpha,\delta}(x)=0$. If $f\in F$ then
$\lim_\beta x_{\alpha,\beta}(f)=v_\alpha(f)$ and
$\lim_\delta\ell_{\alpha,\delta}(f)=v_\alpha(f)-x_\alpha(f)$.
Therefore, by \cite[Theorem V.3.13]{DS58}, for each $\ep>0$ and
$\beta_0$ there is a convex combination of
$\{x_{\alpha,\beta}-x_\alpha-\ell_{\alpha,\delta}\}_{\beta\succ\beta_0}$
satisfying
\[\left\|\sum
a_{\beta,\alpha,\delta}(\beta_0,\ep)(x_{\alpha,\beta}-x_\alpha-\ell_{\alpha,\delta})\right\|<\ep.\]

But then the nets \[\left\{\sum
a_{\beta,\alpha,\delta}(\beta_0,1)x_{\alpha,\beta}\right\}_{\beta_0}\]
are contained in $A$, are uniformly bounded, and
\[w^*-\lim_{\beta_0}\sum
a_{\beta,\alpha,\delta}(\beta_0,1)x_{\alpha,\beta}=x_\alpha.\] We
get a contradiction with the assumption made at the beginning of
the proof.\medskip

We consider the set of all vectors $\{v_\alpha-x_\alpha|_F\}$. It
is clear that it is an unbounded set. We need the following
observation from Convex Geometry.

\begin{lemma}\label{L:UnbConv} Let $\{m(\alpha)\}_{\alpha\in\Omega}\subset\mathbb{R}^n$,
where $\Omega$ is a partially ordered set, be such that
$\limsup_\alpha||m(\alpha)||=\infty$. Then there exist
$0<C<\infty$ and $\alpha'\in\Omega$ such that for each
$\alpha_0\succ\alpha'$ and each $\ep>0$ there is a finitely
non-zero collection $a(\alpha)$ of real numbers supported on
$\alpha\succ\alpha_0$ and satisfying $\sum_\alpha |a(\alpha)|=1$,
$a(\alpha_0)= 1-\ep$, and $\left\|\sum_\alpha
a(\alpha)m(\alpha)\right\|\le C$.
\end{lemma}

We do not specify the norm on $\mathbb{R}^n$ because the lemma
holds for any norm, only the constant $C$ changes.

\begin{proof}[Proof of Lemma \ref{L:UnbConv}] For each $\alpha_0\in\Omega$ consider the closed absolutely convex
hull $M_\alpha$ of $\{m(\alpha)\}_{\alpha\succ\alpha_0}$. By
\cite[Lemma 1.4.2]{Sch93}, each $M_\alpha$ is a (Minkowski) sum of
a compact set $K_\alpha$ and a linear subspace $L_\alpha$.
\medskip

Since $\limsup_\alpha||m(\alpha)||=\infty$, the subspaces
$L_\alpha$ are non-trivial. Also it is clear that
$L_{\alpha_1}\subset L_{\alpha_2}$ for $\alpha_1\succ\alpha_2$.
Since all of these subspaces are finite-dimensional, they
stabilize in the sense that there exists $\alpha'$ such that
$L_{\alpha}= L_{\alpha'}$ for any $\alpha\succeq\alpha'$. Let
$L=L_{\alpha'}(=\cap_\alpha L_\alpha)$. Then $M_\alpha=K_\alpha+L$
for each $\alpha\succeq\alpha'$ and we may assume that
$K_{\alpha}\subset K_{\alpha'}$ (we may assume that all $K_\alpha$
are in the same complement of the subspace $L$, see \cite{Sch93}).
Set $C=\max\{||x||:~x\in K_{\alpha'}\}$.
\medskip

We have $m(\alpha_0)=k(\alpha_0)+\ell(\alpha_0)$, where
$k(\alpha_0)\in K_{\alpha_0}\subset K_{\alpha'}$,
$\ell(\alpha_0)\in L$. Since the vector
$-\frac{1-\ep}{\ep}\ell(\alpha_0)$ is in $L$, it can be
arbitrarily well approximated by absolutely convex combinations of
$\{m(\alpha)\}_{\alpha\succ\alpha_0}$. Therefore there is a
finitely nonzero collection $\{b(\alpha)\}_{\alpha\succ\alpha_0}$
such that $\sum_\alpha|b(\alpha)|=1$ and
\[\left\|\sum_{\alpha\succ\alpha_0}b(\alpha)m(\alpha)+\frac{1-\ep}{\ep}\ell(\alpha_0)\right\|<C.\]
We introduce $a(\alpha)$ by $a(\alpha)=\ep b(\alpha)$ for
$\alpha\succ\alpha_0$, $a(\alpha_0)=1-\ep$ and $a(\alpha)=0$ for
all other $\alpha$. We have
\[\sum_\alpha
a(\alpha)m(\alpha)=(1-\ep)k(\alpha_0)+(1-\ep)\ell(\alpha_0)+ \ep
\sum_{\alpha\succ\alpha_0}b(\alpha)m(\alpha),\] where
$||(1-\ep)k(\alpha_0)||\le (1-\ep)C$ and $||(1-\ep)\ell(\alpha_0)+
\ep \sum_{\alpha\succ\alpha_0}b(\alpha)m(\alpha)||<C\ep$. The
conclusion follows.\end{proof}

We apply Lemma \ref{L:UnbConv} to the set
$\{v_\alpha-x_\alpha|_F\}_\alpha$ and find that there is $C$
(independent of $\alpha$) such that for large enough $\alpha$ and
an arbitrary $\ep>0$ there is a finite combination
\begin{equation}\label{E:combin}
(1-\ep)(v_\alpha-x_\alpha|_F)+\sum_{\delta\succ\alpha}a(\delta)(v_\delta-x_\delta|_F)\end{equation}
having norm $\le C$ and such that $\sum_\delta|a(\delta)|=\ep$.
Using local reflexivity \cite{JRZ71,LR69} we can find a net
$\{p_\gamma\}\subset X^*$ whose weak$^*$ limit is $0$ and whose
restrictions to $F$ converge to the vector \eqref{E:combin}, and
$\sup_\gamma ||p_\gamma||\le C_1$, where $C_1$ does not depend on
$\alpha$.

Then the $(\beta,\gamma)$-net
\begin{equation}
(1-\ep)(x_{\alpha,\beta}-x_\alpha)+\sum_{\delta\succ\alpha}
a(\delta)(x_{\delta,\beta}-x_\delta)-p_\gamma
\end{equation}
is weakly null, where the ordering on  pairs $(\beta,\gamma)$ is
defined as above. Therefore, by \cite[Theorem V.3.13]{DS58}, for
each $\beta_0$ and $\omega>0$ there is a convex combination
satisfying
\begin{equation}\label{E:omega}
\left\|\sum_{\beta\succ\beta_0,\gamma}d_{\alpha,\beta,\gamma,\delta}(\beta_0,\omega)\left((1-\ep)(x_{\alpha,\beta}-x_\alpha)+\sum_{\delta\succ\alpha}
a(\delta)(x_{\delta,\beta}-x_\delta)-p_\gamma\right)\right\|<\omega.
\end{equation}

Consider the net
\[\left\{\sum_{\beta\succ\beta_0,\gamma}d_{\alpha,\beta,\gamma,\delta}(\beta_0,1)\left((1-\ep)x_{\alpha,\beta}+\sum_{\delta\succ\alpha}
a(\delta)x_{\delta,\beta}\right)\right\}_{\beta_0}.\] It is clear
that this net is weak$^*$ convergent to
$(1-\ep)x_\alpha+\sum_{\delta\succ\alpha} a(\delta)x_{\delta}$.
Since $A$ is absolutely convex each element of this net is in $A$.
By \eqref{E:omega}, the elements of this net are norm-bounded
independently of $\alpha$.
\medskip

Now we consider the net
\[\left\{\sum_{\beta\succ\beta_0,\gamma}d_{\alpha,\beta,\gamma,\delta}(\beta_0,1)\left((1-\ep)x_{\alpha,\beta}+\sum_{\delta\succ\alpha}
a(\delta)x_{\delta,\beta}\right)\right\}_{\beta_0,\ep},\] where
$(\beta_1,\ep_1)\succ(\beta_2,\ep_2)$ if and only  if
$\beta_1\succ\beta_2$ and $\ep_1<\ep_2$. It is clear that this net
is weak$^*$ convergent to $x_\alpha$ and its elements are bounded
independently of $\alpha$. This contradicts the assumption made at
the beginning of the proof.
\end{proof}

\end{large}

\begin{thebibliography}{22}
\begin{small}

\bibitem{Alb05} A.\,A.~Albanese, On not open linear continuous operators between Banach spaces,
{\it Note Mat.}, {\bf 25} (2005/06), no. 1, 29--34.

\bibitem{Ban32} S.~Banach, {\it Th\'eorie des op\'erations li\'neaires},
Monografje Matematyczne, Warszawa, 1932.

\bibitem{BDH86} E.~Behrends, S.~Dierolf, P.~Harmand, On a problem of Bellenot and
Dubinsky, {\it Math. Ann.}, {\bf 275} (1986), pp.~337--339.

\bibitem{DJ73}
W.\,J.~Davis, W.\,B.~Johnson, Basic  sequences  and norming
subspaces in  non-quasi-reflexive  Banach  spaces, {\it Israel J.
Math.}, {\bf 14} (1973), 353--367.

\bibitem{DM87}
S.~Dierolf, V.\,B.~Moscatelli, A note on quojections, {\it
Functiones et approximation}, {\bf 17} (1987), 131--138.

\bibitem{Dix48} J.~Dixmier, Sur un th\'eor\`eme de Banach, {\it Duke Math. J.}, {\bf  15}  (1948), 1057--1071.

\bibitem{DS58} N.~Dunford, J.\,T.~Schwartz, {\it Linear
Operators}. Part \textbf{I}: General Theory, New York,
Interscience Publishers, 1958.

\bibitem{GKM10} D.~Garcia, O.\,F.\,K.~Kalenda, M.~Maestre,
Envelopes of open sets and extending holomorphic functions on dual
Banach spaces, {\it J. Math. Anal. Appl.}, {\bf 363} (2010)
663--678, {\tt arXiv:0905.2531}

\bibitem{God78} B.\,V.~Godun, Weak$^*$ derived sets of a set of linear
functionals, {\it Mat. Zametki}, {\bf 23} (1978), 607--616 (in
Russian); English transl. in: {\it Math. Notes.}, {\bf 23} (1978),
333--338.

\bibitem{HMVZ07} P.~Hajek, V.~Montesinos, J.~Vanderwerff, V.~Zizler,
{\it Biorthogonal systems in Banach spaces}, Berlin,
Springer-Verlag, 2007.

\bibitem{JRZ71} W.\,B.~Johnson, H.\,P.~Rosenthal, M.~Zippin, On bases,
finite dimensional decompositions and weaker structures in Banach
spaces, {\it  Israel J. Math.}, {\bf  9} (1971), 488--506.

\bibitem{KL87}
A.\,S.~Kechris, A.~Louveau, {\it Descriptive set theory and the
structure of sets of uniqueness}, Cambridge University Press,
1987.

\bibitem{KKM48} M.\,G.~Krein, M.\,A.~Krasnoselskii, D.\,P.~Milman, On  the  defect
numbers of linear operators  in  a  Banach  space  and  on  some
geometric questions, {\it Sbornik Trudov Inst.  Matem.  AN
Ukrainian SSR}, {\bf 11} (1948), 97--112 (in Russian).

\bibitem{LR69} J.~Lindenstrauss, H.\,P.~Rosenthal, The $\mathcal{L}_{p}$ spaces, {\it Israel J. Math.},
{\bf  7} (1969), 325--349.

\bibitem{LT77} J.~Lindenstrauss, L. Tzafriri, {\it Classical Banach
spaces} {\bf I:} {\it Sequence spaces}, Springer-Ver\-lag, Berlin,
1977.

\bibitem{Lyo88}
R.~Lyons, A new type of sets of uniqueness, {\it Duke Math. J.},
{\bf 57} (1988), 431--458.

\bibitem{Maz30} S.~Mazurkiewicz, Sur la  d\'eriv\'ee  faible  d'un ensemble
de fonctionnelles lin\'eaires, {\it Studia Math.}, {\bf 2} (1930),
68--71.

\bibitem{McG68}
O.\,C.~McGehee, A proof of a statement of Banach about the
weak$^*$ topology, {\it Michigan Math. J.}, {\bf 15} (1968),
135--140.

\bibitem{MM89}
G.~Metafune, V.\,B.~Moscatelli, Quojections and prequojections,
in: {\it Advances in the Theory of Fr\'echet spaces} (ed.:
T.~Terzio\-\v glu),  Kluwer Academic Publishers, Dordrecht,
pp.~235--254, 1989.

\bibitem{MM92} G.~Metafune, V.\,B.~Moscatelli,
Prequojections and their duals, in: {\it Progress in functional
analysis} (eds.: K.\,D.~Bier\-stedt, J.~Bo\-net, J.~Horvath,
M.~Maestre), Pe\~niscola, 1990, North-Holland, Amsterdam,
pp.~215--232, 1992.

\bibitem{Mos87} V.\,B.~Moscatelli, On  strongly  non-norming
subspaces, {\it Note Mat.}, {\bf 7} (1987), 311--314.

\bibitem{Mos90} V.\,B.~Moscatelli,
Strongly    nonnorming    subspaces    and prequojections, {\it
Studia Math.}, {\bf 95} (1990), 249--254.

\bibitem{Ost87}
M.\,I.~Ostrovskii, $w^*$-derived  sets  of  transfinite order of
subspaces of dual Banach spa\-ces, {\it Dokl. Akad. Nauk Ukrain.
SSR}, 1987, no. 10, 9--12 (in Russian and Ukrainian); An English
version of this paper is available at {\tt
http://front.math.ucdavis.edu} as {\tt math/9303203}

\bibitem{Ost98} M.\,I.~Ostrovskii, On prequojections and their duals, {\it Revista Mat.
Univ. Complutense Madrid.}, {\bf 11} (1998), 59--77.

\bibitem{Ost01} M.\,I.~Ostrovskii,
Weak$^*$ sequential closures in Banach space theory and their
applications, in: {\it General Topology in Banach Spaces}, ed. by
T.~Banakh and A.~Plichko, New York, Nova Sci. Publishers, 2001,
pp.~21--34; Available at {\tt http://front.math.ucdavis.edu} as
{\tt math.FA/0203139}

\bibitem{Pel62} A.~Pe\l czy\'nski, A note on the paper of I.~Singer ``Basic
sequences and reflexivity of Banach spaces'', {\it Studia Math.},
{\bf 21} (1961/1962), 371--374.

\bibitem{Pet64} Y.\,I.~Petunin, Conjugate Banach spaces containing subspaces
of zero  characteristic, {\it Dokl.  Akad.  Nauk  SSSR}, {\bf 154}
(1964), 527--529 (in Russian); English transl.: {\it Soviet Math.
Dokl.}, {\bf 5} (1964), 131--133.

\bibitem{PP80} Y.\,I.~Petunin, A.\,N.~Plichko, {\em  The theory of characteristic
of subspaces and its applications}, Vysh\-cha Shko\-la,  Kiev,
1980 (in Russian).

\bibitem{Pia54} I.\,I.~Piatetski-Shapiro,
Supplement to the work ``On the problem of uniqueness of expansion
of a function in a trigonometric series'', {\it Moskov. Gos. Univ.
U\v c. Zap. Mat.}, {\bf 165} (1954), no. 7, 79--97 (in Russian);
English translation in: I.~Piatetski-Shapiro, {\it Selected
works}, Edited by J.~Cogdell, S.~Gindikin, P.~Sarnak, American
Mathematical Society, Providence, RI, 2000.

\bibitem{Pli86}  A.~Plichko, On bounded biorthogonal systems in  some function spaces, {\it Studia
Math.}, {\bf 84} (1986), 25--37.

\bibitem{Pli97}  A. Plichko, Decomposition of Banach space into a direct sum of
separable and reflexive subspaces and Borel maps, {\it Serdica
Math. J.}, {\bf 23} (1997) 335--350.

\bibitem{Raj04} M.~Raja,  Borel properties of linear operators,
{\it J. Math. Anal. Appl.} {\bf 290} (2004), 63--75.

\bibitem{Sai76} J.~Saint-Raymond, Espaces a mod\`ele s\'eparable, {\it Ann. Inst.
Fourier (Grenoble)}, {\bf 26} (1976), 211--256.

\bibitem{Sar68} D.~Sarason, A remark on the weak-star topology of
$\ell^\infty$, {\it Studia Math.}, {\bf 30} (1968), 355--359.

\bibitem{Sch93} R.~Schneider, {\it Convex Bodies: the Brunn--Minkowski Theory}, Encyclopedia of
Mathematics and its Applications, vol. {\bf 44}, Cambridge
University Press, 1993.

\bibitem{Sin62} I.~Singer, Basic sequences and reflexivity of Banach spaces,
{\it Studia Math.}, {\bf 21} (1961/1962), 351--369.

\bibitem{Sin63} I.~Singer, On bases in quasi-reflexive Banach spaces.
{\it Rev. Math. Pures Appl.} (Bucarest) {\bf 8} (1963), 309--311.

\end{small}
\end{thebibliography}
\end{document}